\documentclass[12pt]{article}
\usepackage{bbm}
\usepackage{latexsym}
\usepackage{mathrsfs}
\usepackage{graphicx}
\usepackage{subfigure}
\usepackage{amssymb}
\usepackage{amsmath}
\usepackage{amsthm}
\usepackage{color}
\usepackage{cite}
\usepackage{float}
\usepackage{indentfirst}
\usepackage{anysize}\marginsize{45mm}{45mm}{40mm}{50mm}
\usepackage{caption}
\usepackage{float}
\usepackage{multirow}

\usepackage{enumerate}

\usepackage{latexsym, bm}
\newtheoremstyle{lemma}{\topsep}{\topsep}%
     {}
     {}
     {\bfseries}
     {}
     {0.1em}
     {\thmname{#1}\thmnumber{ #2}\thmnote{ #3}}
\theoremstyle{lemma}  

\newtheorem{theorem}{Theorem}     
\newtheorem{lemma}[theorem]{Lemma}

\newtheorem{conjecture}[theorem]{Conjecture}

\newtheorem{definition}{Definition}

\numberwithin{equation}{section}

\title{ The restricted $h$-connectivity of balanced hypercubes\thanks{This research was partially supported by the National Natural Science Foundation of China (No. 11761056), the Chunhui Project of Ministry of Education (No. Z2017047) and the fundamental research funds for the central universities (No. 2672018ZYGX2018J069)}}

\author{ Huazhong L\"{u}$^{1}$\thanks{Corresponding author.} and Tingzeng Wu$^{2}$\\
{\small $^{1}$School of Mathematical Sciences, University of Electronic Science and Technology of China,} \\
{\small Chengdu, Sichuan 610054, P.R. China}\\
{\small E-mail: lvhz08@lzu.edu.cn}\\
{\small $^{2}$School of Mathematics and Statistics, Qinghai Nationalities University, }\\
{\small Xining, Qinghai 810007, P.R. China} \\
{\small E-mail: mathtzwu@163.com}\\}
\date{}
\begin{document}

\maketitle
\begin{abstract}

The restricted $h$-connectivity of a graph $G$, denoted by $\kappa^h(G)$, is defined as the minimum cardinality of a set of vertices $F$ in $G$, if exists, whose removal disconnects $G$ and the minimum degree of each component of $G-F$ is at least $h$. In this paper, we study the restricted $h$-connectivity of the balanced hypercube $BH_n$ and determine that $\kappa^1(BH_n)=\kappa^2(BH_n)=4n-4$ for $n\geq2$. We also obtain a sharp upper bound of $\kappa^3(BH_n)$ and $\kappa^4(BH_n)$ of $n$-dimension balanced hypercube for $n\geq3$ ($n\neq4$). In particular, we show that $\kappa^3(BH_3)=\kappa^4(BH_3)=12$.

\vskip 0.1 in

\end{abstract}

\section{Introduction}

With the rapid development of VLSI technology, multi-processor systems (MPS) with thousands processors have been available. As the size of MPS increases continuously, the reliability and fault-tolerance become central issues. A MPS is usually modeled by an undirected simple graph with vertices representing processors and edges representing communication links between processors, respectively. Thus, graph parameters can be utilized to evaluate the reliability of a MPS.

The traditional connectivity is an important measure for networks with few processors. A deficiency of the connectivity is the assumption that all the parts of the network can be potentially fail simultaneously. However, in large networks, the possibility that all the vertices incident to a vertex is impossible, indicating their high  resilience. To address the shortcomings of the connectivity stated above, several novel concepts on connectivity has been introduced. Harary \cite{Harary} introduced the conditional connectivity of a connected graph by adding some constraints on the components of the resulting graph after vertex faults occur. After that, the concept of $g$-connectivity was given. On the other hand, Esfahanian and Hakimi \cite{Esfahanian} proposed the restricted connectivity in which all neighbors of any vertex do not fail simultaneously. Later, Latifi et al. \cite{Latifi} generalized it to restricted $h$-connectivity such that every vertex has at least $h$ fault-free neighbors.

A $g$-vertex cut of a graph $G$ is a vertex set $F\subset V(G)$ if $G-F$ is disconnected and each component of $G-F$ contains at least $g+1$ vertices. An $h$-vertex cut of a graph $G$ is a vertex set $F\subset V(G)$ if $G-F$ is disconnected and the minimum degree of each component of $G-F$ is at least $h$. The $g$-connectivity (resp. $h$-connectivity) of $G$, denoted by $\kappa_0^g(G)$ (resp. $\kappa^h(G)$), if exists, is defined as the cardinality of a minimum $g$-vertex cut (resp. $h$-vertex cut). From the definitions above, it is obvious that $\kappa_0^g(G)\leq\kappa_0^{g+1}(G)$ and $\kappa^h(G)\leq \kappa^{h+1}(G)$ if $G$ has $\kappa_0^{g+1}(G)$ and $\kappa^{h+1}(G)$, respectively. In addition, $\kappa_0^0(G)=\kappa^0(G)=\kappa(G)$ for any non-complete graph $G$, and $\kappa_0^1(G)=\kappa^1(G)$. So both of $g$-connectivity and $h$-connectivity are generalizations of the connectivity, which supply more accurate measures to evaluate reliability and fault-tolerance of large networks. Actually, if a network has high $g$-connectivity or $h$-connectivity, then it is more reliable with respect to vertex failure \cite{Esfahanian,Lu2}.

To determine $g$-connectivity and $h$-connectivity of a graph is quite difficult. In fact, no polynomial time algorithm is known for computation of $g$-connectivity and $h$-connectivity for a general graph \cite{Chang,Esfahanian}. The $h$-connectivity \cite{Chen,Li,Ning,Oh,Wan,Ye} and $g$-connectivity \cite{Hsieh,L.Xu,W.Yang,Zhang,Zhu} of given networks are extensively investigated in the literature. For balanced hypercube, Yang \cite{Yang2} showed that $\kappa_0^1(BH_n)=4n-4$ for $n\geq2$. L\"{u} \cite{Lu2} determined that $\kappa_0^2(BH_n)=\kappa_0^3(BH_n)=4n-4$ for $n\geq2$. Recently, Yang et al. \cite{D-Yang} proved that $\kappa_0^4(BH_n)=\kappa_0^5(BH_n)=6n-8$ for $n\geq2$. It seems like more difficult to determine the $h$-connectivity than the $g$-connectivity in $BH_n$, as no research has reported on the $h$-connectivity of $BH_n$ so far. In this paper, we center on the $h$-connectivity of the balanced hypercube.

The balanced hypercube, proposed by Wu and Huang \cite{Wu}, is a novel variant of the well-known hypercube. The balanced hypercube keeps a host of excellent properties of the hypercube, such as bipartite, high symmetry, scalability, etc. It is known that odd-dimension the balanced hypercube has a smaller diameter than that of the hypercube of the same order. Particularly, the balanced hypercube is superior to the hypercube in a sense that it supports an efficient reconfiguration without changing the adjacent relationship among tasks \cite{Wu}. With novel properties mentioned above, it has been attracted considerable attention in literature \cite{Cheng,Hao,Huang2,P.Li,Lu,Lu2,Lu3,Xu,Yang,Yang2,Yang3,Zhou,Zhou2}.

\vskip 0.05 in

The rest of this paper is organized as follows. In Section 2, the definitions of balanced hypercubes and some useful lemmas are presented. The main results of this paper are shown in Section 3. Conclusions are given in Section 4.

\vskip 0.05 in

\section{Preliminaries}

Let $G=(V(G),E(G))$ be a graph, where $V(G)$ is vertex-set of $G$ and $E(G)$ is edge-set of $G$. The number of vertices of $G$ is denoted by $|G|$. The {\em neighborhood} of a vertex $v$ is the set of vertices adjacent to $v$, written as $N_{G}(v)$. Let $F\subseteq V(G)$, we define $N_{G}(F)=\cup_{v\in F}N_{G}(v)-F$. For $A\subset G$, we use $N_{G}(A)$ to denote $N_{G}(V(A))$ briefly. For a subgraph $H\subseteq G$, we use $\delta (H)$ to denote its minimum degree. For other standard graph notations not defined here please refer to \cite{Bondy}.

In what follows, we shall give definitions of the balanced hypercube and some lemmas.
\vskip 0.0 in

\begin{definition}{\bf .}\label{def1}\cite{Wu} An $n$-dimensional balanced hypercube $BH_{n}$ consists of
$2^{2n}$ vertices $(a_{0},\ldots,a_{i-1},a_{i},a_{i+1},\ldots,a_{n-1})$, where $a_{i}\in\{0,1,2,3\}(0\leq
i\leq n-1)$. An arbitrary vertex $v=(a_{0},\ldots,a_{i-1},$
$a_{i},a_{i+1},\ldots,a_{n-1})$ in $BH_{n}$ has the following $2n$ neighbors:

\begin{enumerate}
\item $((a_{0}+1)$ mod $
4,a_{1},\ldots,a_{i-1},a_{i},a_{i+1},\ldots,a_{n-1})$,\\
      $((a_{0}-1)$ mod $ 4,a_{1},\ldots,a_{i-1},a_{i},a_{i+1},\ldots,a_{n-1})$, and
\item $((a_{0}+1)$ mod $ 4,a_{1},\ldots,a_{i-1},(a_{i}+(-1)^{a_{0}})$ mod $
4,a_{i+1},\ldots,a_{n-1})$,\\
      $((a_{0}-1)$ mod $ 4,a_{1},\ldots,a_{i-1},(a_{i}+(-1)^{a_{0}})$ mod $
      4,a_{i+1},\ldots,a_{n-1})$.
\end{enumerate}
\end{definition}

The first coordinate $a_{0}$ of the vertex
$(a_{0},\ldots,a_{i},\ldots,a_{n-1})$ in $BH_{n}$ is defined as {\em inner index}, and
other coordinates $a_{i}$ $(1\leq i\leq n-1)$ {\em outer index}.

\vskip 0.0 in

The recursive structure is one of the most desirable properties in networks. The following definition indicates  recursive property of the balanced hypercube.

\begin{definition}{\bf .}\label{def2}\cite{Wu}
\begin{enumerate}
\item $BH_{1}$ is a $4$-cycle and the vertices are labelled
by $0,1,2,3$ clockwise.
\item $BH_{k+1}$ is constructed from four $BH_{k}$s, which
are labelled by $BH^{0}_{k}$, $BH^{1}_{k}$, $BH^{2}_{k}$,
$BH^{3}_{k}$. For any vertex in $BH_{k}^{i}(0\leq i\leq 3)$, its new labelling in $BH_{k+1}$ is $(a_{0},a_{1},\ldots,a_{k-1},i)$, and it has two new neighbors:
\begin{enumerate}
\item[a)] $BH^{i+1}_{k}:((a_{0}+1)$mod $4,a_{1},\ldots,a_{k-1},(i+1)$mod $4)$ and

$((a_{0}-1)$mod $4,a_{1},\ldots,a_{k-1},(i+1)$mod $4)$ if $a_{0}$ is even.
\item[b)] $BH^{i-1}_{k}:((a_{0}+1)$mod $4,a_{1},\ldots,a_{k-1},(i-1)$mod $4)$ and

$((a_{0}-1)$mod $4,a_{1},\ldots,a_{k-1},(i-1)$mod $4)$ if $a_{0}$ is odd.
\end{enumerate}

\end{enumerate}
\end{definition}

\begin{figure}
\begin{minipage}[t]{0.5\linewidth}
\centering
\includegraphics[width=30mm]{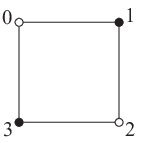}
\caption{$BH_{1}$.} \label{g1}
\end{minipage}
\begin{minipage}[t]{0.5\linewidth}
\centering
\includegraphics[width=60mm]{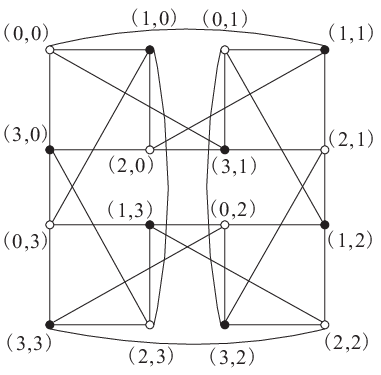}
\caption{$BH_{2}$.} \label{g2}
\end{minipage}
\end{figure}

$BH_{1}$ and $BH_{2}$ are pictured in Figs. \ref{g1} and \ref{g2}, respectively. For brevity, we will omit ``(mod 4)'' in the rest of this paper.

Let $u$ be a neighbor of $v$ in $BH_n$. If $u$ and $v$ differ only from the inner index, then $uv$ is called a $0$-{\em dimension edge}. If $u$ and $v$ differ from $i$th outer index ($1\leq i\leq n-1$), $uv$ is called an $i$-{\em dimension edge}. It can be deduced from Definition \ref{def2} that we can divide $BH_{n}$ into four $BH_{n-1}^{k}$s, $0\leq k \leq 3$, along dimension $n-1$. It is obvious that the edges between $BH_{n-1}^{k}$s are $(n-1)$-dimension edges. Each of $BH_{n-1}^{k}$ is isomorphic to $BH_{n-1}$. For simplicity, we denote $BH_{n-1}^{k}$ by $B^k$ if there is no ambiguity.

The following basic properties of the balanced hypercube will be used in the main result of this paper.

\begin{lemma}\label{bipartite}\cite{Wu}{\bf.}
$BH_{n}$ is bipartite.
\end{lemma}

Since $BH_{n}$ is bipartite, vertices of odd (resp. even) inner index are colored with black (resp. white). The set of all white (resp. black) vertices is denoted by $V_0$ (resp. $V_1$).


\begin{lemma}\label{neighbor}\cite{Wu}{\bf.}
Vertices $u=(a_{0},a_{1},\ldots,a_{n-1})$ and
$v=(a_{0}+2,a_{1},\ldots,a_{n-1})$ in $BH_{n}$ have the same neighborhood.
\end{lemma}

It implies from Lemma \ref{neighbor} that for each vertex $u\in V(BH_n)$, there exists a unique vertex, say $u'$, such that $N(u)=N(u')$.
%
%

\begin{lemma}\label{crossing-edges}{\rm \cite{Yang2}}{\bf.} Assume that $n\geq2$. There exist $4^{n-1}$ edges between $B^i$ and $B^{i+1}$ for each $0\leq i\leq3$.
\end{lemma}

\begin{lemma}\label{g-connectivity-1}\cite{Yang2}{\bf.} $\kappa_0^1(BH_n)=4n-4$ for $n\geq2$.
\end{lemma}

\begin{lemma}\label{g-connectivity-2-3}\cite{Lu2}{\bf.} $\kappa_0^2(BH_n)=\kappa_0^3(BH_n)=4n-4$ for $n\geq2$.
\end{lemma}

\begin{lemma}\label{g-connectivity-4-5}\cite{D-Yang}{\bf.} $\kappa_0^4(BH_n)=\kappa_0^5(BH_n)=6n-8$ for $n\geq2$.
\end{lemma}

\section{Main results}

According to Lemma \ref{g-connectivity-1} and definitions of $g$-connectivity and $h$-connectivity, it is straightforward to obtain the following result.

\begin{lemma}\label{h-connectivity-1}{\bf.} $\kappa^1(BH_n)=4n-4$ for $n\geq2$.
\end{lemma}

The $h$-vertex cut of $BH_n$ for $n\geq2$ is characterized in the following two lemmas.

\begin{lemma}{\bf.}\label{h-vertex-cut-order-even} If $i$-connectivity exists in $BH_n$ for any nonnegative integer $i$ ($2n>i$). Let $F$ be a minimum $i$-vertex cut, then each component of $BH_n-F$ are of even order. Specifically, if $u$ is contained in a component of $BH_n-F$, then $u'$ is contained in the same component.
\end{lemma}

\noindent{\bf Proof.} Let $F$ be a minimum $i$-vertex cut of $BH_n$, then $|F|=\kappa^{i}(BH_n)$. So $BH_n-F$ is disconnected and the minimum degree of each component is no less than $i$. For convenience, assume that $G_1,G_2,\cdots,G_k$ ($k\geq2$) are all components of $BH_n-F$. On the contrary, suppose without loss of generality that there exists a vertex $u'\not\in V(G_1)$ for some $u\in V(G_1)$. Recall that $u$ and $u'$ have the same neighborhood, and $u$ is not connected to the components except $G_1$, then we have $u'\in F$. All neighbors of $u$ are in $G_1$ or $F$, so do $u'$. Thus, $F\setminus\{u'\}$ is a vertex-cut of $BH_n$. Moreover, the subgraph induced by $V(G_1)\cup \{u'\}$, as well as $G_2,\cdots,G_k$ are components of $BH_n-F$ and, clearly, the minimum degree of each component is no less than $i$. Thus, $F\setminus\{u'\}$ is an $i$-vertex cut of $BH_n$, which is a contradiction to the minimality of $F$. It follows that $u'\in V(G_1)$. By the arbitrary choice of $u$ and $G_1$, the result follows immediately.\qed

\begin{lemma}{\bf.}\label{odd-even-equal} If $2i$-connectivity exists for $i\geq1$, then $\kappa^{2i}(BH_n)=\kappa^{2i-1}(BH_n)$.
\end{lemma}

\noindent{\bf Proof.} Obviously, we have $\kappa^{2i-1}(BH_n)\leq\kappa^{2i}(BH_n)$. It remains to show that $\kappa^{2i-1}(BH_n)\geq\kappa^{2i}(BH_n)$. Suppose on the contrary that $\kappa^{2i-1}(BH_n)<\kappa^{2i}(BH_n)$. Let $F$ be a minimum $(2i-1)$-vertex cut of $BH_n$, then $|F|=\kappa^{2i-1}(BH_n)<\kappa^{2i}(BH_n)$. So $BH_n-F$ is disconnected and there exists a component $G_1$ of $BH_n-F$ such that the minimum degree of $G_1$ is $2i-1$. By Lemma \ref{h-vertex-cut-order-even}, $|G_1|$ is even. Moreover, for any vertex $v\in V(G_1)$, if $u\in N_{G_1}(v)$, then $u'\in N_{G_1}(v)$. Thus, by arbitrary choice of $v$ in $G_1$, the degree of each vertex in $G_1$ is even, which is a contradiction. Hence, $\kappa^{2i-1}(BH_n)\geq\kappa^{2i}(BH_n)$. This completes the proof.\qed

Combining Lemma \ref{h-connectivity-1} and Lemma \ref{odd-even-equal}, it is easy to obtain the following theorem.

\begin{theorem}\label{h-connectivity-2}{\bf.} $\kappa^1(BH_n)=\kappa^2(BH_n)=4n-4$ for $n\geq2$.
\end{theorem}

\begin{figure}
\centering
\includegraphics[width=120mm]{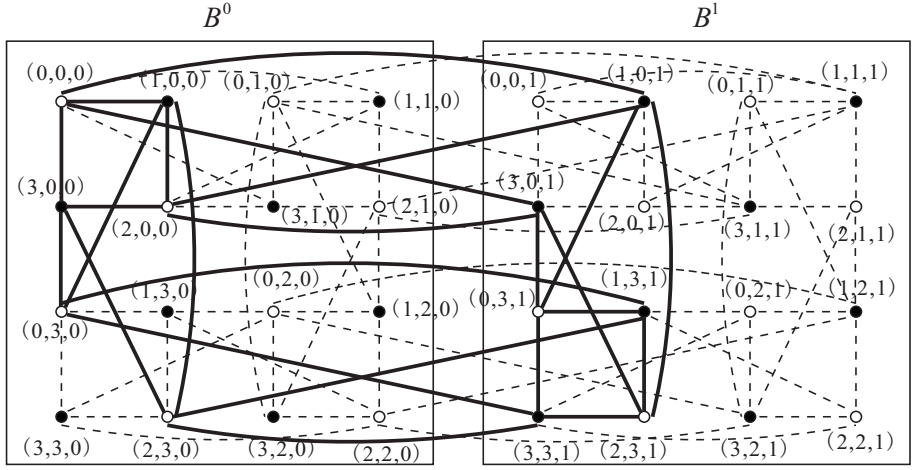}
\caption{A 4-regular subgraph of $BH_{3}$ induced by 12 vertices.} \label{g_counter-example}
\end{figure}

In hypercube and hypercube-like graph, there exist an inequality that, will be utilized to determine their $h$-connectivity \cite{Wu2,Li} as follows. For any subgraph $X\subseteq XQ_n$ and $h\in I_{n}$ satisfying $\delta(X)\geq h$, then $|X|\geq 2^{h}$, where $XQ_n$ means $n$-dimensional hypercube or hypercube-like graph. As a variant of hypercube, one may intuitively speculate that the above inequality holds in $BH_n$. However, it is not true. A subgraph $X$ of $BH_3$ (heavy lines) with $\delta(X)=4$ and $|X|=12$ is shown in Fig. \ref{g_counter-example}, which is obvious a counter-example.

Based on the previous statement, we derive the following lemma.

\begin{lemma}{\bf.}\label{BH3-lower-bound} Let $F\subset V(BH_3)$ with $|F|\leq11$, then either $BH_3-F$ is connected or $BH_3-F$ is disconnected with a component $X$ such that $\delta(X)\leq2$.
\end{lemma}
\noindent{\bf Proof.} It suffices to show that the lemma holds for $|F|=11$. For convenience, let $F^i=B^i\cap F$, $i=0,1,2,3$. Suppose without loss of generality that $|F_0|=\max\{|F_i|:0\leq i\leq3\}$. By Pigeonhole Principle, $|F^0|\geq3$. Obviously, $B^i\cong BH_2$ and $\kappa(B^i)=4$. We consider the following three cases.

\noindent{\bf Case 1.} $|F^0|\geq8$. Since $|F|=11$ and $|F^1|+|F^2|+|F^3|\leq 11-8=3$, each of $B^j-F^j$ is connected, $1\leq j\leq3$. Clearly, $BH_3-B^0-\cup_{j=1}^3 F^j$ is connected.

If $B^0-F^0$ has a single vertex $u$, noting that $u$ has two neighbors $v$ and $w$ in $B^1$ (or $B^3$), furthermore, if $v,w\in F$, then $BH_3-F$ is disconnected with an isolated vertex $u$. If $v\not\in F$ or $w\not\in F$, $BH_3-F$ is connected. So we will only consider the fact that each component of $B^i-F^i$, $0\leq i\leq3$, contains at least two vertices hereafter.

If each component of $B^0-F^0$ contains at least two vertices, then $B^0-F^0$ contains an edge $uv$. Since each vertex in $B^0$ has two neighbors in $B^1$ or $B^3$, $BH_3-F$ is connected.

\noindent{\bf Case 2.} $6\leq|F^0|\leq7$. Then $4\leq|F^1|+|F^2|+|F^3|\leq5$. We further consider the following cases.

\noindent{\bf Case 2.1.} $|F^j|\leq 3$ for each $j\in\{1,2,3\}$. Clearly, each $B^j-F^j$ is connected. Since there are $4^{n-1}$ edges between $B^j$ and $B^{j+1}$ and each vertex of $B^j$ has two neighbors in $B^{j+1}$ or $B^{j-1}$ but not both, there exist at least $4^{n-1}-2|F^j|-2|F^{j+1}|\geq4^{n-1}-2\times3-2\times3>0$ edges between $B^j-F^j$ and $B^{j+1}-F^{j+1}$. So $BH_3-B^0-\cup_{j=1}^3 F^j$ is connected. Let $A$ be an arbitrary component of $B^0-F^0$ (not isolated vertex), then $|N_{BH_3}(A)\setminus V(B^0)|\geq4$. If $|N_{BH_3}(A)\setminus V(B^0)|>4$, then $|N_{BH_3}(A)\setminus V(B^0)|\geq6$, forcing that $A$ is connected to $BH_3-B^0-F$. If $|N_{BH_3}(A)\setminus V(B^0)|=4$ and $|(N_{BH_3}(A)\setminus V(B^0))\cap F|=4$, then $A$ is a component of $BH_3-F$. Note that $|A|\leq4$, it is easy to know that $\delta(A)\leq2$.

\noindent{\bf Case 2.2.} $|F^j|\geq4$ for some $j\in\{1,2,3\}$. If $B^j-F^j$ is connected, the discussion is similar to that of Case 2.1. So we assume that $B^j-F^j$ is disconnected.

If $j=2$, noting $|F^2|\geq4$, $|F^1|\leq1$ and $|F^3|\leq1$. Thus, $BH_3-B^0-F$ is connected. The discussion is also similar to that of Case 2.1.

If $j=1$ or 3, say $j=1$. Without loss of generality, let $A$ be any component of $B^0-F^0$ or $B^1-F^1$, say $B^0-F^0$. then there must exist a vertex $u\in V(A)$ such that $u$ has two neighbors in $B^3-F^3$. Since $|F^3|\leq1$, any component of $B^0-F^0$ is connected to $B^3-F$. By the arbitrariness of $A$, it implies that $BH_3-F$ is connected.

%
%

\noindent{\bf Case 3.} $4\leq|F^0|\leq5$. Then $6\leq|F^1|+|F^2|+|F^3|\leq7$. We further consider the following cases.

\noindent{\bf Case 3.1.} $|F^j|\leq 3$ for each $j\in\{1,2,3\}$. Clearly, each $B^j-F^j$ is connected, thus,  $BH_3-B^0-\cup_{j=1}^3 F^j$ is connected. Let $A$ be any component of $B^0-F^0$ with $|A|\geq2$, then $|N_{BH_3-V(B^0)}(A)|\geq4$.

If $|N_{BH_3-V(B^0)}(A)|=4$. Furthermore, suppose that $|N_{BH_3-V(B^0)}(A)\cap F|=4$, then $A$ is a component of $BH_3-F$. It implies that $|A|\leq 4$, so we have $\delta(A)\leq2$. If $|N_{BH_3-V(B^0)}(A)\cap F|<4$, $A$ is connected to $B^1-F^1$ or $B^3-F^3$, yielding that $BH_3-F$ is connected.

If $|N_{BH_3-V(B^0)}(A)|\geq6$. Since $|F^1|\leq 3$ and $|F^3|\leq 3$, it is obvious that $A$ is connected to $B^1-F^1$ or $B^3-F^3$, thus, $BH_3-F$ is connected.

\noindent{\bf Case 3.2.} $|F^j|\geq 4$ for some $j\in\{1,2,3\}$. If $B^j-F^j$ is connected, the proof is similar to that of Case 3.1. Moreover, if $j=2$, it can be known that $BH_3-B^0-\cup_{j=1}^3 F^j$ is connected, the proof is also similar to that of Case 3.1. Therefore, we may assume that $j=1$ and $B^1-F^1$ is disconnected. It is obvious that $BH_3-V(B^0)-V(B^1)-F^2\cup F^3$ is connected.

If each component of $B^0-F^0$ (resp. $B^1-F^1$) contains a vertex $u$ (resp. $v$) connecting to a neighbor in $BH_3-V(B^0)-V(B^1)-F^2\cup F^3$, then $BH_3-F$ is connected. If not, we may assume that there exists a component $A$ of $B^0-F^0$ such that no vertex of $A$ is connected to $B^3-F^3$. Since $|F^3|\leq3$, $A$ contains exact one black vertex. Additionally, $BH_3$ is bipartite and $B^0$ is 4-regular, so we have $|A|\leq 5$ and $\delta(A)\leq2$. Observe that $4\leq|F^1|\leq5$ and $|F^2|\leq1$, so each component of $B^1-F^1$ contains a vertex $v$ connecting to a neighbor in $B^2-F^2$. If $A$ is connected to a component of $B^1-F^1$, then $BH_3-F$ is connected; if not, $BH_3-F$ is disconnected and contains a component $A$ with $\delta(A)\leq2$.

\noindent{\bf Case 4.} $|F^0|=3$. By proof of Case 2.1, it can be easily known that $BH_3-F$ is connected.

This completes the proof. \qed

\begin{lemma}{\bf.}\label{upper-bound} $\kappa^3(BH_n)=\kappa^4(BH_n)\leq 12n-24$ for $n\geq3$ except $n=4$.
\end{lemma}
\noindent{\bf Proof.} Based on Fig. \ref{g_counter-example}, we derive a subset $T$ of $V(BH_n)$ containing 12 vertices as follows:

$a=(0,0,\cdots,0)$, $a'=(2,0,\cdots,0)$,

$a_1=(1,0,\cdots,0)$, $a_1'=(3,0,\cdots,0)$,

$a_2=(0,3,0,\cdots,0)$, $b'=(2,3,0,\cdots,0)$,

$b=(1,0,1,0,\cdots,0)$, $a_2'=(3,0,1,0\cdots,0)$,

$b_1=(0,3,1,0,\cdots,0)$, $b_2'=(2,3,1,0,\cdots,0)$,

$b_2=(1,3,1,0,\cdots,0)$, $b_1'=(3,3,1,0,\cdots,0)$.

It implies that the graph induced by $T$ is isomorphic to the subgraph of $BH_3$ shown in Fig. \ref{g_counter-example}. For simplicity, denote the subcube $BH_3$ containing $T$ by $H$ in $BH_n$ ($n\geq3$). Obviously, $N_H(T)$ contains the following 12 vertices.

$x=(1,1,0\cdots,0)$, $x'=(3,1,0,\cdots,0)$,

$y=(1,3,0,\cdots,0)$, $y'=(3,3,0,\cdots,0)$,

$z=(0,0,3,0,\cdots,0)$, $z'=(2,0,3,0\cdots,0)$,

$u=(0,0,1,0,\cdots,0)$, $u'=(2,0,1,0,\cdots,0)$,

$v=(0,2,1,0,\cdots,0)$, $v'=(2,2,1,0,\cdots,0)$,

$w=(1,3,2,0,\cdots,0)$, $w'=(3,3,2,0,\cdots,0)$.

For $n\geq4$, $N_{BH_n-H}(T)$ contains the following $12i$ ($1\leq i\leq n-3$) vertices:

$x_i=(1,0,0,\underbrace{0\cdots,0}_{i-1},1,\underbrace{0,\cdots,0}_{n-i-3})$, $x_i'=(3,0,0,\underbrace{0\cdots,0}_{i-1},1,\underbrace{0,\cdots,0}_{n-i-3})$,

$y_i=(1,3,0,\underbrace{0\cdots,0}_{i-1},1,0,\cdots,0)$, $y_i'=(3,3,0,\underbrace{0\cdots,0}_{i-1},1,0,\cdots,0)$,

$z_i=(1,3,1,\underbrace{0\cdots,0}_{i-1},1,0,\cdots,0)$, $z_i'=(3,3,1,\underbrace{0\cdots,0}_{i-1},1,0,\cdots,0)$,

$u_i=(0,0,0,\underbrace{0\cdots,0}_{i-1},3,0,\cdots,0)$, $u_i'=(2,0,0,\underbrace{0\cdots,0}_{i-1},3,0,\cdots,0)$,

$v_i=(0,0,1,\underbrace{0\cdots,0}_{i-1},3,0,\cdots,0)$, $v_i'=(2,0,1,\underbrace{0\cdots,0}_{i-1},3,0,\cdots,0)$,

$w_i=(0,3,1,\underbrace{0\cdots,0}_{i-1},3,0,\cdots,0)$, $w_i'=(2,3,1,\underbrace{0\cdots,0}_{i-1},3,0,\cdots,0)$.

So $|N_{BH_n}(T)|=|N_H(T)|+|N_{BH_n-H}(T)|=12n-24$. It remains to show that $N_{BH_n}(T)$ is a $4$-vertex cut when $n\geq3$. If $n=3$, then $i=0$, it is easy to verify that $N_{BH_3}(T)$ is a $4$-vertex cut. Next we assume that $n\geq4$.

If $n=4$, then $i=1$, implying that $x_1, x_1',y_1, y_1',z_1$ and $z_1'$ are common neighbors of $c_1=(0,3,0,1)$ and $c_1'=(2,3,0,1)$, and $u_1, u_1',v_1, v_1',w_1$ and $w_1'$ are common neighbors of $d_1=(1,0,1,3)$ and $d_1'=(3,0,1,3)$. Therefore, by deleting vertices of $N_H(T)$ and $N_{BH_4-H}(T)$ from $BH_4$, the resulting graph consists of two components $G_1$ and $G_2$, where $G_1$ is the graph induced by $T$. In particular, exact four vertices, i.e. $c_1,c_1',d_1$ and $d_1'$ of $G_2$ have valency two. Hence, the upper bound $12n-24$ is not true for $n=4$.

In what follows, we shall show that the lemma holds for all $n\geq5$. We proceed the proof by induction on $n$.  Firstly, we consider $n=5$. It is obvious that $H\subset B^0$. Specifically, $H$ is a subgraph of $B^0$. Note that, at this time, $B^0=BH_4$, so by above, it can be known that $B^0-T-N_{B^0}(T)$ is connected and exact four vertices, namely $c_1,c_1,d_1$ and $d_1'$, are of valency two in $B^0-T-N_{B^0}(T)$.

By Definition \ref{def2}, $N_{BH_5}(T)\cap V(B^2)=\emptyset$. Clearly, for $i=2$, the vertices in $A=\{x_2, x_2',y_2,y_2',z_2,z_2'\}$ have exact two common neighbors $c_2=(0,3,0,0,1)$ and $c_2'=(2,3,0,0,1)$, and the vertices in $B=\{u_2,u_2',v_2,v_2',w_2,w_2'\}$ have exact two common neighbors $d_2=(1,0,1,0,3)$ and $d_2'=(1,0,1,0,3)$. Similarly, $B^1-A$ (resp. $B^3-B$) is connected with exact two vertices of valency two. As $4^{5-1}-2|A|-2|B|=4^4-2\times12>0$, $BH_5-V(B^0)-N_{BH_5-B^0}(T)$ is connected. In addition, $c_2$ and $c_2'$ (resp. $d_2$ and $d_2'$) have two neighbors in $B^1$ (resp. $B^3$). Thus, the minimum degree of $BH_5-V(B^0)-N_{BH_5-B^0}(T)$ is four. Obviously, $B^0-T-N_{B_0}(T)$ is connected. Additionally, $c_1,c_1',d_1,d_1'\in V(B_0-T-N_{B_0}(T))$, and $c_1$ and $c_1'$ (resp. $d_1$ and $d_1'$) have exact two neighbors in $B^2$. Therefore, the minimum degree of $B^0-T-N_{B_0}(T)$ is four. Hence, $N_{BH_5}(T)$ is a 4-vertex cut of $BH_n$. We assume that the lemma is true for $n-1$ ($n\geq6$). Next we consider $BH_n$.

By induction hypothesis, it follows that $B^0-T-N_{B^0}(T)$ is connected and its minimum degree is no less than 4. Observe that $N_{B^1}(T)=\{x_{n-3},x_{n-3}',y_{n-3},y_{n-3}',z_{n-3}$, $z_{n-3}'\}$ and $N_{B^3}(T)=\{u_{n-3},u_{n-3}',v_{n-3},v_{n-3}',w_{n-3},w_{n-3}'\}$, then $B^1-N_{B^1}(T)$ (resp. $B^3-N_{B^3}(T)$) is also connected and its minimum degree is no less than 4. So $B^1-N_{B^1}(T)$ (resp. $B^3-N_{B^3}(T)$) is connected to $B^2$. Since $4^{n-1}-2|T|-2|N_{BH_n}(T)|+2\times6=4^{n-1}-24-2(12n-24)+12>0$ whenever $n\geq6$, $B^0-T-N_{B^0}(T)$ is connected to $B^1-N_{B^1}(T)$ and $B^3-N_{B^3}(T)$). Therefore, $BH_n-T-N_{BH_n}(T)$ is connected. This completes the proof.\qed

The following statement is straightforward.

\begin{theorem}{\bf.}\label{k-3-4} $\kappa^3(BH_3)=\kappa^4(BH_3)=12$.
\end{theorem}

\section{Conclusions}

In this paper, the more accurate measure of reliability and fault-tolerance, namely $h$-connectivity, of the balanced hypercube is considered. We show that $\kappa^1(BH_n)=\kappa^2(BH_n)=4n-4$ for all $n\geq2$. In addition, we show that $\kappa^{2i-1}(BH_n)=\kappa^{2i}(BH_n)$ for all $1\leq i\leq n-1$. Specifically, we determine  $\kappa^{3}(BH_3)=\kappa^{4}(BH_3)=12$ and derive an upper bound $\kappa^3(BH_n)=\kappa^4(BH_n)\leq 12n-24$ for $n\geq5$. As the upper bound we obtained above is tight for $n=3$, we post the following conjecture as an open problem.

\begin{conjecture}{\bf.} $\kappa^3(BH_n)=\kappa^4(BH_n)=12n-24$ for $n\geq5$.
\end{conjecture}

%

\end{document}